\documentclass{article}
\usepackage{amsmath}

\usepackage{doublespace}

\newtheorem{theorem}{Theorem}

\newtheorem{corollary}[theorem]{Corollary}

\newtheorem{lemma}[theorem]{Lemma}

\newtheorem{proposition}[theorem]{Proposition}
\newtheorem{remark}[theorem]{Remark}

\newenvironment{proof}[1][Proof]{\noindent\textbf{#1.} }{\ \rule{0.5em}{0.5em}}
\input{tcilatex}

\begin{document}

\section{Introduction}

\qquad If $p(x)$ is a polynomial of degree $n\geq 2$ with $n$ distinct real
roots $r_{1}<r_{2}<\cdots <r_{n}$ and critical points $x_{1}<x_{2}<\cdots
<x_{n-1}$, let 
\[
\sigma _{k}=\dfrac{x_{k}-r_{k}}{r_{k+1}-r_{k}},k=1,2,...,n-1. 
\]

$(\sigma _{1},...,\sigma _{n-1})$ is called the \textit{ratio vector} of $p$%
, and $\sigma _{k}$ is called the $k$th ratio. Ratio vectors were first
discussed in \cite{p} and in \cite{a}, where the inequality $\dfrac{1}{n-k+1}%
<\sigma _{k}<\dfrac{k}{k+1},k=1,2,...,n-1$ was derived. For $n=3$ it was
shown in \cite{a} that $\sigma _{1}$ and $\sigma _{2}$ satisfy the
polynomial equation $3(1-\sigma _{1})\sigma _{2}-1=0$. In addition,
necessary and sufficient conditions were given in \cite{H1} for $\left(
\sigma _{1},\sigma _{2}\right) $ to be a ratio vector. For $n=4,$ a
polynomial, $Q$, in three variables was given in \cite{H1} with the property
that $Q\left( \sigma _{1},\sigma _{2},\sigma _{3}\right) =0$ for any ratio
vector $\left( \sigma _{1},\sigma _{2},\sigma _{3}\right) $. It was also
shown that the ratios are monotonic--that is, $\sigma _{1}<\sigma
_{2}<\sigma _{3}$ for any ratio vector $\left( \sigma _{1},\sigma
_{2},\sigma _{3}\right) $. For $n=3,$ $\dfrac{1}{3}<\sigma _{1}<\dfrac{1}{2}%
\,$and $\dfrac{1}{2}<\sigma _{2}<\dfrac{2}{3}$, and thus it follows
immediately that $\sigma _{1}<\sigma _{2}$. The monotonicity of the ratios
does not hold in general for $n\geq 5$(see \cite{H1}). Further results on
ratio vectors for $n=4$ were proved by the author in \cite{H2}. In
particular, necessary and sufficient conditions were given for $\left(
\sigma _{1},\sigma _{2},\sigma _{3}\right) $ to be a ratio vector. We now
want to extend the notion of ratio vector to polynomial like functions of
the form $p(x)=(x-r_{1})^{m_{1}}\cdots (x-r_{N})^{m_{N}}$, where $%
m_{1},...,m_{N}$ are given positive real numbers and $r_{1}<r_{2}<\cdots
<r_{N}$. We extend some of the results and simplify some of the proofs in 
\cite{H1} and in \cite{H2}, and we prove some new results as well. In
particular, we derive more general bounds on the $\sigma _{k}$(Theorem \ref%
{bounds}). Even for $N=3$ or $N=4$, the monotonicity of the ratios does not
hold in general for all positive real numbers $m_{1},...,m_{N}$. We provide
examples below and we also derive necessary and sufficient conditions on $%
m_{1},m_{2},m_{3}$ which imply that $\sigma _{1}<\sigma _{2}$(Theorem \ref%
{T2}).Finally, we prove some general results for any $N$ using Projective
Elimination Theory(see Proposition \ref{P1}). Proposition \ref{P1} can be
used to provide necessary and sufficient conditions for $(\sigma
_{1},...,\sigma _{N-1})$ to be a ratio vector. In particular, we
show(Corollary \ref{C1}) that there is a polynomial $Q\neq 0$ in $N-1$
variables with real coefficients and which does not depend on $%
r_{1},...,r_{N},\sigma _{1},...,\sigma _{N-1}$, such that $Q(\sigma
_{1},...,\sigma _{N-1})=0$ for every ratio vector $(\sigma _{1},...,\sigma
_{N-1})$.

\section{Main Results}

\qquad Throughout, 
\[
p(x)=(x-r_{1})^{m_{1}}\cdots (x-r_{N})^{m_{N}}, 
\]%
where $m_{1},...,m_{N}$ are given positive real numbers with $%
\sum\limits_{k=1}^{N}m_{k}=n$ and $r_{1}<r_{2}<\cdots <r_{N}$. We need the
following lemmas.

\begin{lemma}
\label{L1}$p^{\prime }$ has exactly one root, $x_{k}\in
I_{k}=(r_{k},r_{k+1}),k=1,2,...,N-1$.
\end{lemma}

\begin{proof}
By Rolle's Theorem, $p^{\prime }$ has at least one root in $I_{k}$ for each $%
k=1,2,...,N-1$. Now $\dfrac{p^{\prime }}{p}=\sum\limits_{k=1}^{N}\dfrac{m_{k}%
}{x-r_{k}}$, which has at most $N-1$ real roots since $\left\{ \dfrac{1}{%
x-r_{k}}\right\} _{k=1,...,N}$ is a Chebyshev system.
\end{proof}

Now we define the $N-1$ ratios%
\begin{equation}
\sigma _{k}=\dfrac{x_{k}-r_{k}}{r_{k+1}-r_{k}},k=1,2,...,N-1.  \label{ratio}
\end{equation}

$(\sigma _{1},...,\sigma _{N-1})$ is called the \textit{ratio vector} of $p$%
. We shall derive a system of nonlinear equations in the $\left\{
r_{k}\right\} $ and $\left\{ \sigma _{k}\right\} $. By the product rule, $%
p^{\prime }(x)=(x-r_{1})^{m_{1}-1}\cdots (x-r_{N})^{m_{N}-1}\times $

$\left( \sum\limits_{j=1}^{N}m_{j}\prod\limits_{i=1,i\neq
j}^{N}(x-r_{i})\right) $. Since $p^{\prime }(x)=n(x-r_{1})^{m_{1}-1}\cdots
(x-r_{N})^{m_{N}-1}\times $

$\prod\limits_{k=1}^{N-1}(x-x_{k})$ as well, we have%
\begin{equation}
n\prod\limits_{k=1}^{N-1}(x-x_{k})=\dsum\limits_{j=1}^{N}m_{j}\left(
\prod\limits_{i=1,i\neq j}^{N}(x-r_{i})\right) .  \label{2}
\end{equation}

Let $e_{k}\equiv e_{k}(r_{1},...,r_{n})$ denote the $k$th elementary
symmetric function of the $r_{j}$, $j=1,2,...,n$, starting with $%
e_{0}(r_{1},...,r_{n})=1,e_{1}(r_{1},...,r_{n})=r_{1}+\cdots +r_{n}$, and so
on. Let 
\[
e_{k,j}(r_{1},...,r_{N})=e_{k}(r_{1},...,r_{j-1},r_{j+1},...,r_{N}), 
\]%
that is, $e_{k,j}(r_{1},...,r_{N})$ equals $e_{k}(r_{1},...,r_{N})$ with $%
r_{j}$ removed, $j=1,...,n$. Since $p(x+c)$ and $p(x)$ have the same ratio
vectors for any constant $c$, we may assume that 
\[
r_{2}=0. 
\]

Equating coefficients in (\ref{2}) using the elementary symmetric functions
yields%
\[
ne_{k}(x_{1},...,x_{N-1})=\sum%
\limits_{j=1}^{N}m_{j}e_{k,j}(r_{1},0,r_{3},...,r_{N}),k=1,2,...,N-1. 
\]%
Since $e_{k,j}(r_{1},0,r_{3},...,r_{N})=\left\{ 
\begin{array}{ll}
e_{k,j}(r_{1},r_{3},...,r_{N}) & \text{if }j\neq 2\ \text{and }k\leq N-2 \\ 
e_{k}(r_{1},r_{3},...,r_{N}) & \text{if }j=2 \\ 
0 & \text{if }j\neq 2\text{ and }k=N-1%
\end{array}%
\right. $ we have

\begin{align}
ne_{k}(x_{1},...,x_{N-1})& =m_{2}e_{k}(r_{1},r_{3},...,r_{N})+  \nonumber \\
\sum\limits_{j=1,j\neq 2}^{N}m_{j}e_{k,j}(r_{1},r_{3},...,r_{N}),k&
=1,...,N-2  \label{3} \\
nx_{1}\cdots x_{N-1}& =m_{2}r_{1}r_{3}\cdots r_{N}  \nonumber
\end{align}

Solving (\ref{ratio}) for $x_{k}$ yields 
\begin{equation}
x_{k}=\Delta _{k}\sigma _{k}+r_{k},k=1,2,...,N-1,  \label{xk}
\end{equation}

where $\Delta _{k}=r_{k+1}-r_{k}$. Substituting (\ref{xk}) into (\ref{3})
gives the following equivalent system of equations involving the roots and
the ratios. 
\begin{eqnarray}
ne_{k}((1-\sigma _{1})r_{1},r_{3}\sigma _{2},\Delta _{3}\sigma
_{3}+r_{3},...,\Delta _{N-1}\sigma _{N-1}+r_{N-1}) &=&  \label{4} \\
m_{2}e_{k}(r_{1},r_{3},...,r_{N})+\sum\limits_{j=1,j\neq
2}^{N}m_{j}e_{k,j}(r_{1},r_{3},...,r_{N}),k &=&1,...,N-2  \nonumber \\
n(1-\sigma _{1})r_{1}(r_{3}\sigma _{2})(\Delta _{3}\sigma _{3}+r_{3})\cdots
(\Delta _{N-1}\sigma _{N-1}+r_{N-1}) &=&m_{2}r_{1}r_{3}\cdots r_{N} 
\nonumber
\end{eqnarray}

Note that this system is \textit{homogeneous} in the $r_{k}$ since $\Delta
_{k}\sigma _{k}+r_{k}$ is a linear function of $r_{k}$ for each $k$. This
will be crucial later in our use of \textit{projective} elimination theory.

\subsection{Bounds}

The inequality $\dfrac{1}{n-k+1}<\sigma _{k}<\dfrac{k}{k+1},k=1,2,...,n-1$
was first derived in \cite{p} and later in \cite{a} for polynomials of
degree $n\geq 2$ with $n$ distinct real roots. Critical in proving the
inequality was the root--dragging theorem(see \cite{an}). We now extend this
inequality to the ratios defined in (\ref{ratio}) for functions of the form $%
p(x)=(x-r_{1})^{m_{1}}\cdots (x-r_{N})^{m_{N}}$. First we generalize the
root--dragging theorem. The proof is very similar to the proof in \cite{an}
where $m_{1}=\cdots =m_{N}=1$. For completeness, we provide the details here.

\begin{lemma}
\label{L2}Let $x_{1}<x_{2}<\cdots <x_{N-1}$ be the $N-1$ critical points of $%
p$ lying in $I_{k}=(r_{k},r_{k+1}),k=1,2,...,N-1$. Let $q(x)=(x-r_{1}^{%
\prime })^{m_{1}}\cdots (x-r_{N}^{\prime })^{m_{N}}$, where $r_{k}^{\prime
}>r_{k},k=1,2,...,N-1$ and let $x_{1}^{\prime }<x_{2}^{\prime }<\cdots
<x_{N-1}^{\prime }$ be the $N-1$ critical points of $q$ lying in $%
J_{k}=(r_{k}^{\prime },r_{k+1}^{\prime }),k=1,2,...,N-1$. Then $%
x_{k}^{\prime }>x_{k},k=1,2,...,N-1$.
\end{lemma}

\begin{proof}
Suppose that for some $i,x_{i}^{\prime }<x_{i}$. Now $p^{\prime
}(x_{i})=0\Rightarrow \sum\limits_{k=1}^{N}\dfrac{m_{k}}{x_{i}-r_{k}}=0$ and 
$q^{\prime }(x_{i}^{\prime })=0\Rightarrow \sum\limits_{k=1}^{N}\dfrac{m_{k}%
}{x_{i}^{\prime }-r_{k}^{\prime }}=0$. $r_{k}^{\prime }>r_{k}$ and $%
x_{i}^{\prime }<x_{i}$ implies that%
\begin{equation}
x_{i}^{\prime }-r_{k}^{\prime }<x_{i}-r_{k},k=1,2,...,N-1.  \label{1}
\end{equation}

Since both sides of (\ref{1}) have the same sign, $\dfrac{m_{k}}{%
x_{i}^{\prime }-r_{k}^{\prime }}>\dfrac{m_{k}}{x_{i}-r_{k}},k=1,2,...,N-1$,
which contradicts the fact that $\sum\limits_{k=1}^{N}\dfrac{m_{k}}{%
x_{i}-r_{k}}$ and $\sum\limits_{k=1}^{N}\dfrac{m_{k}}{x_{i}^{\prime
}-r_{k}^{\prime }}$ are both zero.
\end{proof}

\begin{theorem}
\label{bounds}If $\sigma _{1},...,\sigma _{N-1}$ are defined by (\ref{ratio}%
), then 
\begin{equation}
\dfrac{m_{k}}{m_{k}+\cdots +m_{N}}<\sigma _{k}<\dfrac{m_{1}+\cdots +m_{k}}{%
m_{1}+\cdots +m_{k+1}}  \label{5}
\end{equation}
\end{theorem}

\begin{proof}
To obtain an upper bound for $\sigma _{k}$ we use Lemma \ref{L2}. Arguing as
in \cite{a}, we can move the critical point $x_{k}\in (r_{k},r_{k+1})$ as
far to the right as possible by letting $r_{1},...,r_{k-1}\rightarrow r_{k}$
and $r_{k+2},...,r_{N}\rightarrow \infty $. Let $s=m_{1}+\cdots
+m_{k},t=m_{k+2}+\cdots +m_{N}$, and let $%
q_{b}(x)=(x-r_{k})^{s}(x-r_{k+1})^{m_{k+1}}(x-b)^{t}$. Then $q_{b}^{\prime
}(x)=(x-r_{k})^{s}\times $

$\left(
(x-r_{k+1})^{m_{k+1}}t(x-b)^{t-1}+m_{k+1}(x-r_{k+1})^{m_{k+1}-1}(x-b)^{t}%
\right) +$

$s(x-r_{k})^{s-1}(x-r_{k+1})^{m_{k+1}}(x-b)^{t}=$

$(x-r_{k+1})^{m_{k+1}-1}(x-r_{k})^{s-1}(x-b)^{t-1}\times $

$\left( t(x-r_{k+1})(x-r_{k})+m_{k+1}(x-r_{k})(x-b)+s(x-r_{k+1})(x-b)\right) 
$.

$x_{k}$ is the smallest root of the quadratic $t(x-r_{k+1})(x-r_{k})+$

$m_{k+1}(x-r_{k})(x-b)+s(x-r_{k+1})(x-b)=\left( m_{k+1}+t+s\right) x^{2}+$

$\left( -tr_{k+1}-tr_{k}-m_{k+1}r_{k}-m_{k+1}b-sr_{k+1}-sb\right)
x+tr_{k+1}r_{k}+sr_{k+1}b+m_{k+1}r_{k}b$. As $b\rightarrow \infty ,x_{k}$
increases and approaches the root of $(-m_{k+1}-s)x+sr_{k+1}+m_{k+1}r_{k}$.
Thus $x_{k}\uparrow \dfrac{sr_{k+1}+m_{k+1}r_{k}}{m_{k+1}+s}\Rightarrow $

$\sigma _{k}\uparrow \left( \dfrac{sr_{k+1}+m_{k+1}r_{k}}{m_{k+1}+s}%
-r_{k}\right) /(r_{k+1}-r_{k})=$

$\dfrac{sr_{k+1}+m_{k+1}r_{k}-r_{k}(m_{k+1}+s)}{(m_{k+1}+s)(r_{k+1}-r_{k})}=%
\dfrac{s}{m_{k+1}+s}=\dfrac{m_{1}+\cdots +m_{k}}{m_{1}+\cdots +m_{k+1}}$.
Similarly, to obtain a lower bound for $\sigma _{k}$, move the critical
point $x_{k}\in (r_{k},r_{k+1})$ as far to the left as possible by letting $%
r_{k+2},...,r_{N}\rightarrow r_{k+1}$ and $r_{1},...,r_{k-1}\rightarrow
-\infty $. By considering $%
q_{b}(x)=(x-r_{k})^{m_{k}}(x-r_{k+1})^{s}(x+b)^{t} $, where $%
s=m_{k+1}+\cdots +m_{N}$ and $t=m_{1}+\cdots +m_{k-1}$, one obtains $\sigma
_{k}\downarrow \dfrac{m_{k}}{m_{k}+\cdots +m_{N}}$.
\end{proof}

\subsection{N = 3}

The following Theorem generalizes (\cite{H1},Theorem 1). Throughout, 
\[
n=m_{1}+m_{2}+m_{3} 
\]

\begin{theorem}
\label{T1}Let $p(x)=(x-r_{1})^{m_{1}}(x-r_{2})^{m_{2}}(x-r_{3})^{m_{3}}$.
Then $\left( \sigma _{1},\sigma _{2}\right) $ is a ratio vector if and only
if $\dfrac{m_{1}}{n}<\sigma _{1}<\dfrac{m_{1}}{m_{1}+m_{2}},\dfrac{m_{2}}{%
m_{2}+m_{3}}<\sigma _{2}<\dfrac{m_{1}+m_{2}}{n}$, and $\sigma _{2}=\dfrac{%
m_{2}}{n(1-\sigma _{1})}$
\end{theorem}

\begin{proof}
To prove the necessity part, 
\begin{equation}
\dfrac{m_{1}}{n}<\sigma _{1}<\dfrac{m_{1}}{m_{1}+m_{2}},\dfrac{m_{2}}{%
m_{2}+m_{3}}<\sigma _{2}<\dfrac{m_{1}+m_{2}}{n}  \label{6}
\end{equation}%
follows from Theorem \ref{bounds} with $N=3$. With $N=3$ (\ref{4}) becomes 
\begin{eqnarray*}
n(r_{2}\sigma _{1}+(r_{3}-r_{2})\sigma _{2}+r_{2})
&=&m_{1}(r_{2}+r_{3})+m_{2}r_{3}+m_{3}r_{2} \\
nr_{2}\sigma _{1}((r_{3}-r_{2})\sigma _{2}+r_{2}) &=&m_{1}(r_{2}r_{3})
\end{eqnarray*}%
Since $p(cx)$ and $p(x)$ have the same ratios when $c>0$, in addition to $%
r_{1}=0$ we may also assume that $r_{2}=1$. Let $r_{2}=1$ and $r_{3}=r$ to
obtain

\begin{eqnarray}
\left( n\sigma _{2}-m_{1}-m_{2}\right) r+n(\sigma _{1}-\sigma _{2})+m_{2}
&=&0  \notag \\
\left( n\sigma _{1}\sigma _{2}-m_{1}\right) r+\allowbreak n\sigma _{1}\left(
1-\sigma _{2}\right) &=&0  \label{7}
\end{eqnarray}

Note that $n\sigma _{2}-m_{1}-m_{2}\neq 0$ since $\sigma _{2}<\dfrac{%
m_{1}+m_{2}}{n}$ by (\ref{6}) and $n\sigma _{1}\sigma _{2}-m_{1}\neq 0$
since $\sigma _{1}\sigma _{2}<\dfrac{m_{1}}{m_{1}+m_{2}}\dfrac{m_{1}+m_{2}}{n%
}=\dfrac{m_{1}}{n}$ by (\ref{6}). Hence we can solve each equation in (\ref%
{7}) for $r$ to obtain $r=\dfrac{n(\sigma _{2}-\sigma _{1})-m_{2}}{n\sigma
_{2}-m_{1}-m_{2}}$ and $r=n\sigma _{1}\dfrac{\sigma _{2}-1}{n\sigma
_{1}\sigma _{2}-m_{1}}$. Equating these two expressions yields $\dfrac{%
n(\sigma _{2}-\sigma _{1})-m_{2}}{n\sigma _{2}-m_{1}-m_{2}}=n\sigma _{1}%
\dfrac{\sigma _{2}-1}{n\sigma _{1}\sigma _{2}-m_{1}}\Rightarrow \left(
n\sigma _{1}\sigma _{2}-n\sigma _{2}+m_{2}\right) \left( m_{1}-n\sigma
_{1}\right) =0\Rightarrow \sigma _{2}=\dfrac{m_{2}}{n(1-\sigma _{1})}$ since 
$m_{1}-n\sigma _{1}\neq 0$. To prove sufficiency, suppose that $u$ is any
real number with $\dfrac{m_{1}}{n}<u<\dfrac{m_{1}}{m_{1}+m_{2}}$. We want to
show that $\left( u,v\right) $ is a ratio vector, where $v=\dfrac{m_{2}}{%
n(1-u)}$. Let $r=u\dfrac{(1-u)n-m_{2}}{m_{1}-(m_{1}+m_{2})u}$. $u<\dfrac{%
m_{1}}{m_{1}+m_{2}}\Rightarrow m_{1}-(m_{1}+m_{2})u>0$ and $(1-u)n-m_{2}>%
\dfrac{m_{2}}{m_{1}+m_{2}}n-m_{2}=\allowbreak m_{2}\dfrac{n-m_{1}-m_{2}}{%
m_{1}+m_{2}}=\dfrac{m_{2}m_{3}}{m_{1}+m_{2}}>0$. Also, $u\dfrac{(1-u)n-m_{2}%
}{m_{1}-(m_{1}+m_{2})u}>1\iff u\left( (1-u)n-m_{2}\right)
>m_{1}-(m_{1}+m_{2})u\iff u\left( (1-u)n-m_{2}\right) -\left(
m_{1}-(m_{1}+m_{2})u\right) =\allowbreak \left( 1-u\right) \left(
un-m_{1}\right) >0$, which holds since $\dfrac{m_{1}}{n}<u$ and $u<\dfrac{%
m_{1}}{m_{1}+m_{2}}<1$. Thus $r>1$. Now let $%
p(x)=x^{m_{1}}(x-1)^{m_{2}}(x-r)^{m_{3}}$. Then the ratios of $p,\sigma _{1}$
and $\sigma _{2}$, must satisfy (\ref{7}) by (\ref{4}) with $N=3$. Thus $%
r=\sigma _{1}\dfrac{(1-\sigma _{1})n-m_{2}}{m_{1}-(m_{1}+m_{2})\sigma _{1}}$%
. For any $r$, (\ref{7}) with $\sigma _{1}=u$ and $\sigma _{2}=v$ is
equivalent to $f(u)=0$, where $f(u)=-nu^{2}+\left(
rm_{1}+rm_{2}+n-m_{2}\right) u-rm_{1}$. Then $f\left( \dfrac{m_{1}}{n}%
\right) =\allowbreak m_{1}\left( r-1\right) \dfrac{-m_{3}}{n}<0$ and $%
f\left( \dfrac{m_{1}}{m_{1}+m_{2}}\right) =\allowbreak \dfrac{m_{1}m_{2}m_{3}%
}{\left( m_{1}+m_{2}\right) ^{2}}>0$. If $r=u\dfrac{(1-u)n-m_{2}}{%
m_{1}-(m_{1}+m_{2})u}$, then $f(u)=\allowbreak 0$. Since $%
\lim\limits_{u\rightarrow \pm \infty }f(u)=-\infty $, $f$ \ has exactly one
solution $\dfrac{m_{1}}{n}<u<\dfrac{m_{1}}{m_{1}+m_{2}}$. Now $\dfrac{m_{1}}{%
n}<\sigma _{1}<\dfrac{m_{1}}{m_{1}+m_{2}}$ and $f(\sigma _{1})=0$ as well.
Thus $u=\sigma _{1}$. That finishes the proof of Theorem \ref{T1}.
\end{proof}

\subsubsection{Monotonicity}

For $m_{1}=m_{2}=m_{3}=1$, Theorem \ref{bounds} yields $\dfrac{1}{3}<\sigma
_{1}<\dfrac{1}{2}\,$and $\dfrac{1}{2}<\sigma _{2}<\dfrac{2}{3}$, and thus it
follows immediately that $\sigma _{1}<\sigma _{2}$. $\sigma _{1}<\sigma _{2}$
does not hold in general for all positive real numbers(or even positive
integers) $m_{1},m_{2},$ and $m_{3}$. For example, if $m_{1}=6$, $m_{2}=1$,
and $m_{3}=2,$ then it is not hard to show that $\sigma _{2}<\sigma _{1}$
for all $r_{1}<r_{2}<r_{3}$. Also, if $m_{1}=4$, $m_{2}=3$, and $m_{3}=6$,
then $\sigma _{1}<\sigma _{2}$ for certain $r_{1}<r_{2}<r_{3}$, while $%
\sigma _{2}<\sigma _{1}$ for other $r_{1}<r_{2}<r_{3}$. For $%
p(x)=x^{4}(x-1)^{3}\left( x+\dfrac{1}{2}-\dfrac{1}{2}\sqrt{13}\right)
^{6},\sigma _{1}=\sigma _{2}=\dfrac{1}{2}-\dfrac{1}{26}\sqrt{13}$. One can
easily derive sufficient conditions on $m_{1},m_{2},m_{3}$ which imply that $%
\sigma _{1}<\sigma _{2}$ for all $r_{1}<r_{2}<r_{3}$. For example, if $%
m_{1}m_{3}<m_{2}^{2}$, then $\dfrac{m_{1}}{m_{1}+m_{2}}<\dfrac{m_{2}}{%
m_{2}+m_{3}}$, which implies that $\sigma _{1}<\sigma _{2}$ by (\ref{6}).
Also, if $m_{1}+m_{3}<3m_{2}$, then $n<4m_{2}$, which implies that $\sigma
_{2}=\dfrac{m_{2}}{n(1-\sigma _{1})}>\dfrac{1}{4(1-\sigma _{1})}\geq \sigma
_{1}$ since $4x(1-x)\leq 1$ for all real $x$. We shall now derive necessary
and sufficient conditions on $m_{1},m_{2},m_{3}$ which imply that $\sigma
_{1}<\sigma _{2}$.

\begin{theorem}
\label{T2}$\sigma _{1}<\sigma _{2}$ for all $r_{1}<r_{2}<r_{3}$ if and only
if

(A) $m_{2}^{2}+m_{1}(m_{2}-m_{3})>0$ and

(B) $m_{2}\geq \dfrac{2m_{1}m_{3}}{n}$ or

(C) $\dfrac{n}{4}<m_{2}<\dfrac{2m_{1}m_{3}}{n}$
\end{theorem}

\begin{proof}
As noted above, we may assume that $%
p(x)=x^{m_{1}}(x-1)^{m_{2}}(x-r)^{m_{3}},r>1$. A simple computation shows
that 
\begin{equation*}
\sigma _{1}=\dfrac{1}{2n}\left( (n-m_{3})r-n-m_{2}-\sqrt{A}\right) +1
\end{equation*}%
\begin{equation*}
\sigma _{2}=\dfrac{\dfrac{1}{2n}\left( (n-m_{3})r-n-m_{2}+\sqrt{A}\right) }{%
r-1}
\end{equation*}%
where 
\begin{equation*}
A=(m_{1}+m_{2})^{2}r^{2}+2(m_{2}m_{3}-m_{1}n)r+(m_{1}+m_{3})^{2}
\end{equation*}%
Thus $\sigma _{2}-\sigma _{1}=\dfrac{\dfrac{1}{2n}\left( (n-m_{3})r-n-m_{2}+%
\sqrt{A}\right) }{r-1}-$

$\dfrac{1}{2n}\left( (n-m_{3})r-n-m_{2}-\sqrt{A}\right) -1=\allowbreak $

$\dfrac{1}{2}\dfrac{-\left( n-m_{3}\right) r^{2}-\left(
-n+2m_{3}-m_{2}\right) r-2m_{2}+\sqrt{A}r}{n\left( r-1\right) }>0$ when $%
r>1\iff \sqrt{A}r>\allowbreak \left( n-m_{3}\right) r^{2}+\left(
-n+2m_{3}-m_{2}\right) r+2m_{2}\iff $

$Ar^{2}>\left( \allowbreak \left( n-m_{3}\right) r^{2}+\left(
-n+2m_{3}-m_{2}\right) r+2m_{2}\right) ^{2}\iff $

$4\left( r-1\right) \left( \allowbreak \left(
m_{2}^{2}+m_{1}m_{2}-m_{1}m_{3}\right) r^{2}+\left(
m_{2}m_{3}-m_{1}m_{2}-m_{2}^{2}\right) \allowbreak r+m_{2}^{2}\right) >0\iff
h(r)>0$, where 
\begin{equation*}
h(r)=\allowbreak \left( m_{2}^{2}+m_{1}(m_{2}-m_{3})\right)
r^{2}+m_{2}\left( m_{3}-m_{2}-m_{1}\right) \allowbreak r+m_{2}^{2}.
\end{equation*}%
We want to determine necessary and sufficient conditions on $%
m_{1},m_{2},m_{3}$ which imply that $h(r)>0$ for all $r>1$. First (A) is
clearly a necessary condition, so we assume that (A) holds. Let $r_{0}=-%
\dfrac{1}{2}m_{2}\dfrac{m_{3}-m_{2}-m_{1}}{m_{2}^{2}+m_{1}m_{2}-m_{1}m_{3}}$
be the unique root of $h^{\prime }$. It suffices to determine when $%
r_{0}\leq 1$ or when $h\left( r_{0}\right) >0$. Now $r_{0}\leq 1\iff 2\left(
m_{2}^{2}+m_{1}m_{2}-m_{1}m_{3}\right) \geq $

$m_{2}\allowbreak \left( m_{2}+m_{1}-m_{3}\right) \iff \allowbreak
m_{2}^{2}+m_{1}m_{2}+m_{2}m_{3}-2m_{1}m_{3}\geq 0$, which yields (B). If $%
m_{2}^{2}+m_{1}m_{2}+m_{2}m_{3}-2m_{1}m_{3}<0$, then it is necessary and
sufficient that $h\left( r_{0}\right) =\dfrac{1}{4}m_{2}^{2}\left(
m_{1}+m_{2}+m_{3}\right) \dfrac{m_{1}+m_{3}-3m_{2}}{%
-m_{2}^{2}-m_{1}m_{2}+m_{1}m_{3}}>0\iff m_{1}+m_{3}-3m_{2}>0$. That yields
(C).
\end{proof}

As noted above, if $m_{1}=m_{2}=m_{3}=1$, then $\sigma _{1}<\sigma _{2}$.
The following corollary is a slight generalization of that and follows
immediately from Theorem \ref{T2}.

\begin{corollary}
Suppose that $m_{1}=m_{2}=m_{3}=m>0$. Then $\sigma _{1}<\sigma _{2}$ for all 
$r_{1}<r_{2}<r_{3}$.
\end{corollary}

\subsection{N = 4}

Throughout, 
\[
n=m_{1}+m_{2}+m_{3}+m_{4} 
\]

To simplify the notation, we use $\sigma _{1}=u,\sigma _{2}=v,$and $\sigma
_{3}=w$ for the ratios. For $N=4$ Theorem \ref{bounds} yields

\begin{eqnarray}
\dfrac{m_{1}}{n} &<&u<\dfrac{m_{1}}{m_{1}+m_{2}}  \nonumber \\
\dfrac{m_{2}}{m_{2}+m_{3}+m_{4}} &<&v<\dfrac{m_{1}+m_{2}}{m_{1}+m_{2}+m_{3}}
\label{12} \\
\dfrac{m_{3}}{m_{3}+m_{4}} &<&w<\dfrac{m_{1}+m_{2}+m_{3}}{n},  \nonumber
\end{eqnarray}

In \cite{H2} necessary and sufficient conditions were given for $\left(
\sigma _{1},\sigma _{2},\sigma _{3}\right) $ to be a ratio vector when $%
m_{1}=m_{2}=m_{3}=1$. We now give a simpler proof than that given in \cite%
{H2} and which also generalizes to general positive real numbers $%
m_{1},m_{2},$ and $m_{3}$. The proof here for $N=4$ does not require the use
of Grobner bases as in \cite{H2}, though we shall use Grobner bases later in
this paper to prove some results for $N$ in general.

\begin{theorem}
\label{T3}Let 
\[
D\equiv D(u,v,w)=\left\vert 
\begin{array}{ll}
n(w-v)-m_{3} & n(1-w)-m_{4} \\ 
n\left( u-1\right) v\left( 1-w\right) & n\left( u-1\right) vw+m_{2}%
\end{array}%
\right\vert , 
\]%
$D_{1}\equiv D_{1}(u,v,w)=\allowbreak \left( nu-m_{1}\right) \left(
m_{2}-nvw\left( 1-u\right) \right) $, $D_{2}\equiv D_{2}(u,v,w)=\left(
nu-m_{1}\right) nv\left( 1-u\right) \left( 1-w\right) $, and 
\begin{gather*}
R\equiv R(u,v,w)= \\
\tfrac{nv(1-w)D_{1}^{2}+\left( nvw-m_{1}-m_{2}\right) D_{1}D_{2}+\left(
n\left( 1-u\right) \left( w-v-1\right) +m_{2}+m_{4}\right) D_{1}D+\left(
nw(u-1)+m_{2}+m_{3}\right) D_{2}D}{\left( nu-m_{1}\right) \left(
m_{2}-nv(1-u)\right) },
\end{gather*}%
which is a polynomial in $u,v,$ and $w$ of degree $7$. Then $(u,v,w)\in \Re
^{3}$ is a ratio vector of $%
p(x)=(x-r_{1})^{m_{1}}(x-r_{2})^{m_{2}}(x-r_{3})^{m_{3}}(x-r_{4})^{m_{4}}$
if and only if $0<D_{1}(u,v,w)<D_{2}(u,v,w),D(u,v,w)>0$, and $%
R(u,v,w)=\allowbreak 0$.
\end{theorem}

\begin{proof}
$(\Longleftarrow $ Suppose first that $(u,v,w)$ is a ratio vector of $%
p(x)=(x-r_{1})^{m_{1}}(x-r_{2})^{m_{2}}(x-r_{3})^{m_{3}}(x-r_{4})^{m_{4}}$.
Since $p(x+c)$ and $p(x)$ have the same ratio vectors for any constant $c$,
we may assume that $r_{2}=0$, and thus the equations (\ref{4}) hold with $%
N=4 $. In addition, since $p(cx)$ and $p(x)$ have the same ratio vectors for
any constant $c>0$, we may also assume that $r_{1}=-1$. In addition, we let $%
r_{3}=r$ and $r_{4}=s$, so that $0<r<s$. Then (\ref{4}) becomes 
\begin{gather}
\left( n(w-v)-m_{3}\right) \allowbreak r+(n(1-w)-m_{4})s=nu-m_{1}  \label{8}
\\
nv(1-w)r^{2}+\left( nvw-m_{1}-m_{2}\right) rs+\left( n\left( 1-u\right)
\left( w-v-1\right) +m_{2}+m_{4}\right) \allowbreak r+  \label{9} \\
\left( nw(u-1)+m_{2}+m_{3}\right) s=0  \notag \\
nv\left( u-1\right) \left( 1-w\right) r+\left( nvw\left( u-1\right)
+m_{2}\right) s=0  \label{10}
\end{gather}

In particular, (\ref{8})--(\ref{10}) must be consistent. Eliminating $r$ and 
$s$ from (\ref{8}) and (\ref{10}) yields

$\left( nv\left( u-1\right) \left( 1-w\right) (n(1-w)-m_{4})-\left(
n(w-v)-m_{3}\right) \left( nvw\left( u-1\right) +m_{2}\right) \right) s=$

$\left( nu-m_{1}\right) nv\left( u-1\right) \left( 1-w\right) $ or $%
Ds=\left( nu-m_{1}\right) nv\left( 1-u\right) \left( 1-w\right) $. Note that 
$nu-m_{1}>0,1-u>0,v>0,$ and $1-w>0$ by \ref{12}. Thus $D\neq 0$ and by
Cramer's Rule, 
\begin{equation}
r=\dfrac{D_{1}(u,v,w)}{D(u,v,w)},s=\dfrac{D_{2}(u,v,w)}{D(u,v,w)},
\label{11}
\end{equation}

where $D_{1}(u,v,w)=\left\vert 
\begin{array}{ll}
nu-m_{1} & n(1-w)-m_{4} \\ 
0 & nvw\left( u-1\right) +m_{2}%
\end{array}%
\right\vert =\allowbreak $

$\left( nu-m_{1}\right) \left( m_{2}-nvw\left( 1-u\right) \right) $, and $%
D_{2}(u,v,w)=$

$\left\vert 
\begin{array}{ll}
n(w-v)-m_{3} & nu-m_{1} \\ 
nv\left( u-1\right) \left( 1-w\right) & 0%
\end{array}%
\right\vert =\allowbreak $

$\left( nu-m_{1}\right) nv\left( 1-u\right) \left( 1-w\right) $. $D_{2}>0$
and $s>0$ implies that $D>0$, which in turn implies that $D_{1}>0$ since $%
r>0 $. $r<s\Rightarrow D_{1}<D_{2}$. Now substitute the expressions for $r$
and $s$ in (\ref{11}) into (\ref{9}). Clearing denominators gives 
\begin{gather}
nv(1-w)\left( D_{1}(u,v,w)\right) ^{2}+\left( nvw-m_{1}-m_{2}\right)
D_{1}(u,v,w)D_{2}(u,v,w)+  \label{13} \\
\left( n\left( 1-u\right) \left( w-v-1\right) +m_{2}+m_{4}\right)
D_{1}(u,v,w)D(u,v,w)+  \notag \\
\left( nw(u-1)+m_{2}+m_{3}\right) D_{2}(u,v,w)D(u,v,w)=0.  \notag
\end{gather}%
Factoring the LHS of (\ref{13}) yields $\left( nu-m_{1}\right) \left(
nv(1-u)-m_{2}\right) R(u,v,w)=0$. Also, (\ref{10}) and $r<s$ implies that $%
\dfrac{m_{2}}{n}-vw\left( 1-u\right) <v\left( 1-u\right) \left( 1-w\right)
\Rightarrow \dfrac{m_{2}}{n}<vw\left( 1-u\right) +v\left( 1-u\right) \left(
1-w\right) =\allowbreak v\left( 1-u\right) \Rightarrow $%
\begin{equation}
v(1-u)>\dfrac{m_{2}}{n}.  \label{14}
\end{equation}%
Thus $m_{2}-nv(1-u)\neq 0$, which implies that $R(u,v,w)=0$.

$(\Longrightarrow $ Now suppose that $u,v,$ and $w$ are real numbers with $%
0<D_{1}(u,v,w)<D_{2}(u,v,w),D(u,v,w)>0$, and $R(u,v,w)=\allowbreak 0$. Let $%
r=\dfrac{D_{1}(u,v,w)}{D(u,v,w)}$ and $s=\dfrac{D_{2}(u,v,w)}{D(u,v,w)}$.
Then $0<r<s$ and it follows as above that $(r,s,u,v,w)$ satisfies (\ref{8}%
)--(\ref{10}). Let $x_{1}=u-1,x_{2}=rv,$ and $x_{3}=(s-r)w+r$. Then (\ref{3}%
) must hold since (\ref{3}) and (\ref{4}) are an equivalent system of
equations. Let $p(x)=(x+1)^{m_{1}}x^{m_{2}}(x-r)^{m_{3}}(x-s)^{m_{4}}$.
Working backwards, it is easy to see that (\ref{2}) must hold and hence $%
x_{1},x_{2},$ and $x_{3}$ must be the critical points of $p$. Since $u=%
\dfrac{x_{1}-(-1)}{0-(-1)}$, $v=\dfrac{x_{2}-0}{r-0}$, and $w=\dfrac{x_{3}-r%
}{s-r}$, $(u,v,w)$ is a ratio vector of $p$.
\end{proof}

\begin{remark}
As noted in \cite{H2} for the case when $m_{1}=m_{2}=m_{3}=m_{4}=1$, the
proof above shows that if $(u,v,w)$ is a ratio vector, then there are 
\textit{unique} real numbers $0<r<s$ such that the polynomial $%
p(x)=(x+1)^{m_{1}}x^{m_{2}}(x-r)^{m_{3}}(x-s)^{m_{4}}$ has $(u,v,w)$ as a
ratio vector. For general $N$ we have the following.
\end{remark}

\textbf{Conjecture:} Let $p(x)=(x+1)^{m_{1}}x^{m_{2}}(x-r_{3})^{m_{3}}\cdots
(x-r_{N})^{m_{N}},$

$q(x)=(x+1)^{m_{1}}x^{m_{2}}(x-s_{3})^{m_{3}}\cdots (x-s_{N})^{m_{N}}$,
where $0<r_{3}<\cdots <r_{N}$ and $0<s_{3}<\cdots <s_{N}$. Suppose that $p$
and $q$ have the same ratio vectors. Then $p=q$.

As with $N=3$, it was shown in \cite{H1} that $m_{1}=m_{2}=m_{3}=m_{4}=1%
\Rightarrow \sigma _{1}<\sigma _{2}<\sigma _{3}$. Not suprisingly, this does
not hold for general positive real numbers $m_{1},m_{2},m_{3},$ and $m_{4}$.
Indeed it is possible that $\sigma _{1}>\sigma _{3}$. For example, if $%
p(x)=(x+1)^{3/2}x(x-4)^{\sqrt{2}}(x-6)^{2}$, then $\sigma _{1}>\sigma
_{3}>\sigma _{2}$.

\begin{theorem}
\label{T4}Suppose that $m_{1}+m_{4}\leq \min \left\{
3m_{2}-m_{3},3m_{3}-m_{2}\right\} $. Then $\sigma _{1}<\sigma _{2}<\sigma
_{3}$.
\end{theorem}

\begin{proof}
$m_{1}+m_{4}\leq 3m_{2}-m_{3}\Rightarrow n\leq 4m_{2}$. By (\ref{14}) in the
proof of Theorem \ref{T3}, $v(1-u)>\dfrac{1}{4}$. Thus $\dfrac{v}{u}>\dfrac{1%
}{4u(1-u)}\geq 1$ since $u(1-u)\leq 1$. By letting $%
r_{1}=r<r_{2}=-1<r_{3}=0<r_{4}=s$ one can derive equations similar to (\ref%
{3}) with $N=4$. The third equation becomes

$\left( m_{3}-nw\left( 1-u\right) \left( 1-v\right) \right) r+nwu\left(
1-v\right) =0\Rightarrow $

$r=\dfrac{nwu(1-v)}{nw\left( 1-v\right) \left( 1-u\right) -m_{3}}$. $%
r<-1\Rightarrow \dfrac{1}{r}>-1\Rightarrow \dfrac{nw\left( 1-v\right) \left(
1-u\right) -m_{3}}{nwu(1-v)}>-1\Rightarrow nw\left( 1-v\right) \left(
1-u\right) -m_{3}>-nwu(1-v)\Rightarrow $

$nw\left( 1-v\right) \left( 1-u\right) +nwu(1-v)>m_{3}\Rightarrow nw\left(
1-v\right) >m_{3}\Rightarrow \dfrac{w}{v}>\dfrac{m_{3}}{nv(1-v)}$. Now $%
m_{1}+m_{4}\leq 3m_{3}-m_{2}\Rightarrow n\leq 4m_{3}$. Thus $\dfrac{w}{v}>%
\dfrac{1}{4v(1-v)}\geq 1$.
\end{proof}

As with $N=3$, we have the following generalization of the case when $%
m_{1}=m_{2}=m_{3}=m_{4}=1$, which follows immediately from Theorem \ref{T4}

\begin{corollary}
Suppose that $m_{1}=m_{2}=m_{3}=m_{4}=m>0$. Then $\sigma _{1}<\sigma
_{2}<\sigma _{3}$.
\end{corollary}

We do not derive necessary and sufficient conditions in general on $%
m_{1},m_{2},m_{3},m_{4}$\ which imply that $\sigma _{1}<\sigma _{2}<\sigma
_{3}$.

\subsection{Results for General \ N}

\textit{\ }We note again that throughout $m_{1},...,m_{N}$ are given
positive real numbers.

\begin{lemma}
\label{L3}If $\sigma _{1}=\cdots =\sigma _{N-1}=1$, then the only solution
of (\ref{4}) is $r_{1}=r_{3}=\cdots =r_{N}=0.$
\end{lemma}

\begin{proof}
Let $\sigma _{1}=\cdots =\sigma _{N-1}=1$ in (\ref{xk}), which gives $%
x_{k}=r_{k+1},k=1,...,N-1$. Note that we have assumed that $r_{2}=0$ in
deriving (\ref{4}). Since (\ref{4}) is equivalent to (\ref{3}), we can
substitute $x_{k}=r_{k+1}$ into the last equation in (\ref{3}), which yields 
$nr_{2}\cdots r_{N}=m_{2}r_{1}r_{3}\cdots r_{N}$. Since $r_{2}=0$, the
latter equation implies that $r_{1}r_{3}\cdots r_{N}=0$. Since $r_{l}=0$ for
some $l,l\neq 2$,we may reorder the roots, if necessary, so that $r_{1}=0$.
If $N=2$ we are finished. Otherwise, follow the steps as above with $%
p(x)=x^{m_{1}+m_{2}}(x-r_{3})^{m_{3}}\cdots (x-r_{N})^{m_{N}}$. That is,
replace $m_{1}$ with $m_{1}+m_{2}$, $r_{k}$ with $r_{k+1}$, and $N$ with $%
N-1 $. We would then obtain $r_{3}\cdots r_{N}=0$. Reordering the roots
again, if necessary, $r_{3}=0$. After a finite number of steps, $%
r_{1}=r_{3}=\cdots =r_{N}=0$.
\end{proof}

\begin{proposition}
\label{P1}There are nonzero polynomials in $N-1$ variables, $Q_{1},...,Q_{l}$%
, which do not depend on $r_{1},r_{3},...,r_{N},\sigma _{1},...,\sigma
_{N-1} $, with the following property. There are complex numbers $%
r_{1},r_{3},...,r_{N}$ such that

$(r_{1},r_{3},...,r_{N},\sigma _{1},...,\sigma _{N-1})$ is a solution of (%
\ref{4}) if and only if

$Q_{1}(\sigma _{1},...,\sigma _{N-1})=\cdots =Q_{l}(\sigma _{1},...,\sigma
_{N-1})=0$.
\end{proposition}

\begin{proof}
For each $k=1,2,...,N-2$, let%
\begin{gather*}
f_{k}(r_{1},r_{3},...,r_{N},\sigma _{1},...,\sigma
_{N-1})=m_{2}e_{k}(r_{1},r_{3},...,r_{N})+ \\
\sum\limits_{j=1,j\neq 2}^{N}m_{j}e_{k,j}(r_{1},r_{3},...,r_{N})- \\
ne_{k}((1-\sigma _{1})r_{1},r_{3}\sigma _{2},\Delta _{3}\sigma
_{3}+r_{3},...,\Delta _{N-1}\sigma _{N-1}+r_{N-1}), \\
k=1,2,...,N-2,f_{N-1}(r_{1},r_{3},...,r_{N},\sigma _{1},...,\sigma _{N-1})=
\\
m_{2}r_{1}r_{3}\cdots r_{N}-n(1-\sigma _{1})r_{1}(r_{3}\sigma _{2})(\Delta
_{3}\sigma _{3}+r_{3})\cdots (\Delta _{N-1}\sigma _{N-1}+r_{N-1}),
\end{gather*}

and 
\begin{gather*}
V_{a}=V(f_{1},...,f_{N-1})= \\
\left\{ 
\begin{array}{c}
(r_{1},r_{3},...,r_{N},\sigma _{1},...,\sigma _{N-1})\in
C^{2n-2}:f_{k}(r_{1},r_{3},...,r_{N},\sigma _{1},...,\sigma _{N-1})=0, \\ 
k=1,...,N-1%
\end{array}%
\right\} .
\end{gather*}%
Then the solutions in $C^{2n-2}$ of the system (\ref{4}) are precisely the
points of the \textit{affine} variety $V_{a}$. It is much more useful,
however, to view $r_{k},k=2,...,N$ as variables in projective $N-2$ space, $%
P^{N-2}$. Define the \textit{projective ideal }%
\begin{equation*}
I=\langle f_{1},...,f_{N-1}\rangle ,
\end{equation*}%
the ideal generated by $f_{1},...,f_{N-1}$ in $P^{N-2}\times C^{N-1}$ and 
\begin{gather*}
V=V(f_{1},...,f_{N-1})= \\
\left\{ (r_{1},r_{3},...,r_{N},\sigma _{1},...,\sigma _{N-1})\in
P^{N-2}\times C^{N-1}:f_{k}=0,k=1,...,N-1\right\} .
\end{gather*}%
Now we want to use Projective Elimination Theory(see \cite{clo}, Chapter 8).
Define the \textit{projective elimination ideal } 
\begin{equation*}
\hat{I}=\left\{ f\in C[\sigma _{1},...,\sigma _{N-1}]:\text{for each }j\text{%
, there is }e_{j}\geq 0\text{ with }r_{j}^{e_{j}}f\in I\right\} .
\end{equation*}%
Let $\pi :P^{N-2}\times C^{N-1}\rightarrow C^{N-1}$ be the projection map.
Since each of the polynomials $f_{1},...,f_{N-1}$ is \textit{homogeneous} in 
$r_{1},r_{3},...,r_{N}$, by the Projective Extension Theorem(\cite{clo},
page 389, Theorem 6), 
\begin{equation}
\pi (V)=V(\hat{I})  \label{piv}
\end{equation}%
Since $V(\hat{I})$ is an affine variety(by definition) contained in $C^{N-1}$%
, by (\ref{piv}) $\pi (V)$ is also an affine variety contained in $C^{N-1}$.
By Lemma \ref{L3}, $\pi (V)$ cannot be all of $C^{N-1}$ since $%
(0,...,0)\notin P^{N-2}$, which implies that $\pi (V)$ is a \textit{proper}
affine variety. That finishes the proof since $(\sigma _{1},...,\sigma
_{N-1})\in \pi (V)$ if and only if there are complex numbers $%
r_{1},r_{3},...,r_{N}$ such that $(r_{1},r_{3},...,r_{N},\sigma
_{1},...,\sigma _{N-1})$ is a solution of (\ref{4}).
\end{proof}

\begin{corollary}
\label{C1}There is a polynomial $Q\neq 0$ in $N-1$ variables with real
coefficients and which does not depend on $r_{1},...,r_{N},\sigma
_{1},...,\sigma _{N-1}$, such that $Q(\sigma _{1},...,\sigma _{N-1})=0$ for
every ratio vector $(\sigma _{1},...,\sigma _{N-1}).$
\end{corollary}

\begin{proof}
If $(\sigma _{1},...,\sigma _{N-1})$ is a ratio vector of $%
p(x)=(x-r_{1})^{m_{1}}\cdots (x-r_{N})^{m_{N}},r_{1}<r_{2}<\cdots <r_{N}$,
then we may assume that $r_{2}=0$. Then there are complex numbers $%
r_{1},r_{3},...,r_{N}$ such that $(r_{1},r_{3},...,r_{N},\sigma
_{1},...,\sigma _{N-1})$ is a solution of (\ref{4}). By Proposition \ref{P1}%
, $Q(\sigma _{1},...,\sigma _{N-1})=0$ for some polynomial $Q\neq 0$ in $N-1$
variables which does not depend on $r_{1},r_{3},...,r_{N},\sigma
_{1},...,\sigma _{N-1}$. \ Since the ratios of a polynomial with real roots
must be real, by taking real and imaginary parts of $Q(\sigma
_{1},...,\sigma _{n-1})$, it follows immediately that one can assume that $Q$
has \textit{real} coefficients.
\end{proof}

\begin{remark}
With a little more effort, one can show that $Q$ has \textit{integer}
coefficients.
\end{remark}

\begin{remark}
Corollary \ref{C1} can be proven without using Grobner bases or Projective
Elimination Theory. Instead one can use some theory and facts about the
Krull dimension of an ideal. However, this approach is really not much
shorter and more importantly, it does not yield the sufficiency part of
Proposition \ref{P1}, which is a stronger result than Corollary \ref{C1}.
Proposition \ref{P1} can be used to obtain sufficient conditions for $%
(\sigma _{1},...,\sigma _{N-1})$ to be a ratio vector. Additional
restrictions on $(\sigma _{1},...,\sigma _{N-1})$ are needed which would
force $r_{1},...,r_{n}$ to be real and distinct, with $r_{1}<\cdots <r_{n}$.
One can also check if a particular $(\sigma _{1},...,\sigma _{N-1})$
satisfies $Q(\sigma _{1},...,\sigma _{N-1})=0$. One then knows, without
solving (\ref{4}), that there are complex numbers $r_{1},r_{3},...,r_{N}$
such that $(r_{1},r_{3},...,r_{N},\sigma _{1},...,\sigma _{N-1})$ is an
exact solution of (\ref{4}). Then (\ref{4}) can be solved numerically to see
if $r_{1},r_{3},...,r_{N}$ are real with $r_{1}<r_{3}<\cdots <r_{N}$.
\end{remark}

\end{document}